\magnification=\magstep1
\input amstex
\documentstyle{amsppt}
\loadeusm

\input xy
\xyoption{all}
\newdir{ >}{{}*!/-10pt/@{>}}

\topmatter
\title Two-vector bundles and forms of elliptic cohomology \endtitle
\author Nils A. Baas, Bj{\o}rn Ian Dundas and John Rognes \endauthor 
\dedicatory Dedicated to Graeme Segal on the occasion of his 60th
birthday \enddedicatory
\address Department of Mathematical Sciences, NTNU, Trondheim, Norway
\endaddress
\email baas\@math.ntnu.no \endemail
\address Department of Mathematical Sciences, NTNU, Trondheim, Norway
\endaddress
\email dundas\@math.ntnu.no \endemail
\address Department of Mathematics, University of Oslo, Norway \endaddress
\email rognes\@math.uio.no \endemail
\endtopmatter

\define\Ar{\operatorname{Ar}}
\define\A{\Bbb A}
\define\C{\Bbb C}
\define\E{\eusm E}
\define\Fun{\operatorname{Fun}}
\define\F{\eusm F}
\define\Gr{\operatorname{Gr}}
\define\Iso{\operatorname{Iso}}
\define\I{\eusm I}
\define\Map{\operatorname{Map}}
\define\M{\eusm M}
\define\N{\Bbb N}
\define\Q{\Bbb Q}
\define\R{\Bbb R}
\define\SS{\Bbb S}
\define\Spec{\operatorname{Spec}}
\define\Spt{\operatorname{Spt}}
\define\TMF{\operatorname{TMF}}
\define\Tel{\operatorname{Tel}}
\define\T{\eusm T}
\define\U{\eusm U}
\define\V{\eusm V}
\define\W{\eusm W}
\define\Z{\Bbb Z}
\define\colim{\operatornamewithlimits{colim}}
\define\hofib{\operatorname{hofib}}
\define\telecom{\operatorname{telecom}}
\define\tmf{\operatorname{tmf}}
\define\trc{\operatorname{trc}}
\define\twoVect{\operatorname{2-Vect}}
\redefine\B{\eusm B}

\hyphenation{co-ho-mo-logy}

\document

\head 1. Introduction \endhead

The work to be presented in this paper has been inspired by several
of Professor Graeme Segal's papers.  Our search for a geometrically
defined elliptic cohomology theory with associated elliptic objects
obviously stems from his Bourbaki seminar \cite{Se88}.  Our readiness
to form group completions of symmetric monoidal categories by passage to
algebraic $K$-theory spectra derives from his Topology paper \cite{Se74}.
Our inclination to invoke 2-functors to the 2-category of 2-vector spaces
generalizes his model for topological $K$-theory in terms of functors
from a path category to the category of vector spaces.  We offer him
our admiration.

\medskip

Among all generalized (co-)homology theories, a few hold a special
position because they are, in some sense, geometrically defined.
For example, de\,Rham cohomology of manifolds is defined in terms of
cohomology classes of closed differential forms, topological $K$-theory of
finite CW complexes is defined in terms of equivalence classes of complex
vector bundles, and complex bordism is defined in terms of bordism classes
of maps from stably complex manifolds.  The geometric origin of these
theories makes them particularly well suited to the analysis of many
key problems.  For example, Chern--Weil theory associates differential
forms related to the curvature tensor to manifolds with a connection,
whose de\,Rham cohomology classes are the Chern classes of the tangent
bundle of the manifold.  The Atiyah--Segal index theory \cite{AS68}
associates formal differences of vector bundles to parametrized
families of Fredholm operators, arising e.g.~from complexes of elliptic
pseudo-differential operators, and their isomorphism classes live in
topological $K$-theory.  Moduli spaces of isomorphism classes of solutions
to e.g.~Yang--Mills gauge-theoretic problems can generically yield maps
from suitably structured manifolds, with well-defined bordism classes in
the corresponding form of bordism homology.  On the other hand, Quillen's
theorem that the coefficient ring $\pi_*(MU)$ for complex bordism theory
is the Lazard ring that corepresents (commutative $1$-dimensional)
formal group laws has no direct manifold-geometric interpretation,
and may seem to be a fortuitous coincidence in this context.

From the chromatic point of view of stable homotopy theory, related
to the various periodicity operators $v_n$ for $n\ge0$ that act in
many cohomology theories, these three geometrically defined cohomology
theories detect an increasing amount of information.  De\,Rham cohomology
or real singular cohomology sees only rational phenomena, because for each
prime $p$ multiplication by $p = v_0$ acts invertibly on $H^*(X; \R)$.
Topological $K$-theory only picks up Bott periodic phenomena, because
multiplication by the Bott class $u \in \pi_2(KU)$ acts invertibly on
$KU^*(X)$, and $u^{p-1} = v_1$ for each prime $p$.  Complex bordism
$MU_*(X)$ instead detects all levels of periodic phenomena.  We can
say that real cohomology, topological $K$-theory and complex bordism
have {\it chromatic filtration}~$0$, $1$ and~$\infty$, respectively.
A precise interpretation of this is that the spectra $H\R$ and $KU$
are Bousfield local with respect to the Johnson--Wilson spectra $E(n)$
for $n=0$ and $1$, respectively, while $MU$ is not $E(n)$-local for any
finite~$n$.  Traditionally, an elliptic cohomology theory is a complex
oriented Landweber exact cohomology theory associated to the formal
group law of an elliptic curve.  It will have chromatic filtration~$2$
when the elliptic curve admits a supersingular specialization, and so any
cohomology theory of chromatic filtration~$2$ might loosely be called a
form of elliptic cohomology.  However, the formal group law origin of
traditional elliptic cohomology is not of a directly geometric nature, and so
there has been some lasting interest in finding a truly geometrically
defined form of elliptic cohomology.

\medskip

It is the aim of the present paper to introduce a geometrically
defined cohomology theory that is essentially of chromatic
filtration~$2$, or more precisely, a connective form of such a theory.
It therefore extends the above list of distinguished cohomology theories
one step beyond topological $K$-theory, to a theory that will detect
$v_2$-periodic phenomena, but will ignore the complexity of all higher
$v_n$-periodicities for $n\ge3$.

\medskip

The theory that we will present is represented by the algebraic $K$-theory
spectrum $K(\V)$ of the Kapranov--Voevodsky 2-category of 2-vector spaces
\cite{KV94}.  A 2-vector space is much like a complex vector space, but
with all occurrences of complex numbers, sums, products and equalities
replaced by finite-dimensional complex vector spaces, direct sums, tensor
products and coherent isomorphisms, respectively.  It is geometrically
defined in the sense that the $0$-th cohomology group $K(\V)^0(X)$
of a space $X$ can be defined in terms of equivalence classes of {\it
2-vector bundles\/} over~$X$ (or more precisely, over the total space $Y$
of a Serre fibration $Y \to X$ with acyclic homotopy fibers, i.e., an
acyclic fibration).  Cf.~theorem~4.10.  A 2-vector bundle over~$X$ is a
suitable bundle of categories, defined much like a complex vector bundle
over~$X$, but subject to the same replacements as above.  The previously
studied notion of a gerbe over~$X$ with band $\C^*$ is a special case of
a 2-vector bundle, corresponding in the same way to a complex line bundle.

We conjecture in~5.1 that the spectrum $K(\V)$ is equivalent to the
algebraic $K$-theory spectrum $K(ku)$ of the connective topological
$K$-theory spectrum $ku$, considered as a ``brave new ring'',
i.e., as an $\SS$-algebra.  This is a special case of a more general
conjecture, where for a symmetric bimonoidal category $\B$ (which is a
generalization of a commutative semi-ring) we compare the category of
finitely generated free modules over~$\B$ to the category of finitely
generated free modules over the commutative $\SS$-algebra $\A = \Spt(\B)$
(which is a generalization of a commutative ring) associated to $\B$.
The conjecture amounts to a form of ``positive thinking'', asserting
that for the purpose of forming algebraic $K$-theory spectra it should
not matter whether we start with a semi-ring-like object (such as the
symmetric bimonoidal category $\B$) or the ring-like object given by
its additive Grothendieck group completion (such as the commutative
$\SS$-algebra~$\A$).  This idea originated with Marcel B{\"o}kstedt, and
we are indebted to him for suggesting this approach.  We have verified
the conjecture in the case of actual commutative semi-rings, interpreted
as symmetric bimonoidal categories that only have identity morphisms,
and view this as strong support in favor of the conjecture.

Continuing, we know that $K(ku)$, or rather a spectrum very closely
related to it, is essentially of chromatic filtration~$2$.  For connective
spectra, such as all those arising from algebraic $K$-theory, there is
a more appropriate and flexible variation of the chromatic filtration
that we call the {\it telescopic complexity\/} of the spectrum;
cf.~definition~6.1.  For example, integral and real cohomology have
telescopic complexity~$0$, connective and periodic topological $K$-theory
have telescopic complexity~$1$, and traditional elliptic cohomology has
telescopic complexity~$2$.

It is known, by direct nontrivial calculations \cite{AR02}, that
$K(\ell^\wedge_p)$ has telescopic complexity~$2$, where $\ell^\wedge_p$
is the connective $p$-complete Adams summand of topological $K$-theory
and $p\ge5$.  The use of the Adams summand in place of the full
connective $p$-complete topological $K$-theory spectrum $ku^\wedge_p$,
as well as the hypothesis $p\ge5$, are mostly technical assumptions that
make the calculations manageable, and it seems very likely that also
$K(ku^\wedge_p)$ will have telescopic complexity~$2$ for any prime $p$.
It then follows from \cite{D97}, if we assume the highly respectable
Lichtenbaum--Quillen conjecture for $K(\Z)$ at $p$, that also $K(ku)$
has telescopic complexity~$2$.  In this sense we shall allow ourselves
to think of $K(ku)$, and conjecturally $K(\V)$, as a connective form of
elliptic cohomology.

The definition of a 2-vector bundle is sufficiently explicit that it may
carry independent interest.  In particular, it may admit notions of {\it
connective structure\/} and {\it curving}, generalizing the notions for
gerbes \cite{Br93, \S5.3}, such that to each 2-vector bundle $\E$ over~$X$
with connective structure there is an associated virtual vector bundle $H$
over the free loop space $\eusm L X = \Map(S^1, X)$, generalizing the {\it
anomaly line bundle\/} for gerbes \cite{Br93, \S6.2}.  If $\E$ is equipped
with a curving, there probably arises an {\it action functional\/} for
oriented compact surfaces over $X$ (loc.~cit.), providing a construction
of an {\it elliptic object\/} over~$X$ in the sense of Segal \cite{Se88}.
Thus 2-vector bundles over~$X$ (with extra structure) may have naturally
associated elliptic objects over~$X$.  However, we have not yet developed
this theory properly, and shall therefore postpone its discussion to
a later joint paper, which will also contain proofs of the results
announced in the present paper.  Some of the basic ideas presented here
were sketched by the first author in \cite{Ba98}.

\medskip

The paper is organized as follows.  In \S2 we define a charted
2-vector bundle of rank~$n$ over a space $X$ with respect to an open
cover $\U$ that is indexed by a suitably partially ordered set $\I$.
This corresponds to a Steenrod-style definition of a fiber bundle, with
standard fiber the category $\V^n$ of $n$-tuples of finite-dimensional
complex vector spaces, chosen trivializations over the chart domains in
$\U$, gluing data that compare the trivializations over the intersection
of two chart domains and coherence isomorphisms that systematically
relate the two possible comparisons that result over the intersection
of three chart domains.  We also discuss when two such charted 2-vector
bundles are to be viewed as equivalent, i.e., when they define the same
abstract object.

In \S3 we think of a symmetric bimonoidal category $\B$ as a generalized
semi-ring, and make sense of the algebraic $K$-theory $K(\B)$ of its
2-category of finitely generated free ``modules'' $\B^n$.  We define the
weak equivalences $\B^n \to \B^n$ to be given by a monoidal category
$\M = GL_n(\B)$ of {\it weakly invertible\/} matrices over~$\B$,
cf.~definition~3.6, in line with analogous constructions for simplicial
rings and $\SS$-algebras \cite{W78}.  It is a key point that we allow
$GL_n(\B)$ to contain more matrices than the strictly invertible ones,
of which there are too few to yield an interesting theory.  We also
present an explicit bar construction $B\M$ that is appropriate for such
monoidal categories.  Our principal example is the symmetric bimonoidal
category $\V$ of finite-dimensional complex vector spaces under direct
sum and tensor product, for which the modules $\V^n$ are the 2-vector
spaces of Kapranov and Voevodsky.

In \S4 we bring these two developments together, by showing that the
equivalence classes of charted 2-vector bundles of rank~$n$ over a
(reasonable) space $X$ is in natural bijection (theorem~4.5) with
the homotopy classes of maps from $X$ to the geometric realization
$|BGL_n(\V)|$ of the bar construction on the monoidal category of weakly
invertible $n \times n$ matrices over~$\V$.  The group of homotopy classes
of maps from $X$ to the algebraic $K$-theory space $K(\V)$ is naturally
isomorphic (theorem~4.10) to the Grothendieck group completion of the
abelian monoid of virtual 2-vector bundles over~$X$, i.e., the 2-vector
bundles $\E \downarrow Y$ over spaces~$Y$ that come equipped with an
acyclic fibration $a \: Y \to X$.  Hence the contravariant homotopy
functor represented by $K(\V)$ is geometrically defined, in the sense
that virtual 2-vector bundles over~$X$ are the (effective) cycles for
this functor at $X$.

In \S5 we compare the algebraic $K$-theory of the generalized semi-ring
$\B$ to the algebraic $K$-theory of its additive group completion.
To make sense of the latter as a ring object, as is necessary to form its
algebraic $K$-theory, we pass to structured ring spectra, i.e., to the
commutative $\SS$-algebra $\A = \Spt(\B)$.  We propose that the resulting
algebraic $K$-theory spectra $K(\B)$ and $K(\A)$ are weakly equivalent
(conjecture 5.1), and support this assertion by confirming that it holds
true in the special case of a discrete symmetric bimonoidal category $\B$,
i.e., a commutative semi-ring in the usual sense.  In the special case
of 2-vector spaces the conjecture asserts that $K(\V)$ is the algebraic
$K$-theory $K(ku)$ of connective topological $K$-theory $ku$ viewed as
a commutative $\SS$-algebra.

In \S6 we relate the spectrum $K(ku)$ to the algebraic $K$-theory
spectrum $K(\ell^\wedge_p)$ of the connective $p$-complete Adams summand
$\ell^\wedge_p$ of $ku^\wedge_p$.  The latter theory $K(\ell^\wedge_p)$
is known (theorem~6.4, \cite{AR02}) to have telescopic complexity~$2$,
and this section makes it plausible that also the former theory $K(ku)$
has telescopic complexity~$2$, and hence is a connective form of
elliptic cohomology.  Together with conjecture 5.1 this says that (a)
the generalized cohomology theory represented by $K(ku)$ is geometrically
defined, because its $0$-th cohomology group, which is then represented
by $K(\V)$, is defined in terms of formal differences of virtual 2-vector
bundles, and (b) that it has telescopic complexity~$2$, meaning that
it captures one more layer of chromatic complexity than topological
$K$-theory does.

\head 2. Charted two-vector bundles \endhead

\definition{Definition 2.1}
Let $X$ be a topological space.  An {\it ordered open cover} $(\U, \I)$
of $X$ is a collection $\U = \{U_\alpha \mid \alpha \in \I\}$ of open
subsets $U_\alpha \subset X$, indexed by a partially ordered set $\I$,
such that
\roster
\item
the $U_\alpha$ cover $X$ in the sense that $\bigcup_\alpha U_\alpha =
X$, and
\item
the partial ordering on $\I$ restricts to a total ordering on each finite
subset $\{\alpha_0, \dots, \alpha_p\}$ of $\I$ for which the intersection
$U_{\alpha_0 \dots \alpha_p} = U_{\alpha_0} \cap \dots \cap U_{\alpha_p}$
is nonempty.
\endroster
The partial ordering on $\I$ makes the nerve of the open cover $\U$ an
ordered simplicial complex, rather than just a simplicial complex.
We say that $\U$ is a {\it good cover\/} if each finite intersection
$U_{\alpha_0 \dots \alpha_p}$ is either empty or contractible.
\enddefinition

\definition{Definition 2.2}
Let $X$ be a topological space, with an ordered open cover $(\U, \I)$, and
let $n \in \N = \{0, 1, 2, \dots\}$ be a non-negative integer.
A {\it charted 2-vector bundle} $\E$ of rank~$n$ over~$X$ consists of
\roster
\item
an $n \times n$ matrix
$$
E^{\alpha\beta} = (E^{\alpha\beta}_{ij})_{i,j=1}^n
$$
of complex vector bundles over~$U_{\alpha\beta}$, for each pair $\alpha
< \beta$ in $\I$, such that over each point $x \in U_{\alpha\beta}$
the integer matrix of fiber dimensions
$$
\dim(E^{\alpha\beta}_x) = (\dim E^{\alpha\beta}_{ij,x})_{i,j=1}^n
$$
is invertible, i.e., has determinant $\pm1$, and
\item
an $n \times n$ matrix
$$
\phi^{\alpha\beta\gamma} = (\phi^{\alpha\beta\gamma}_{ik})_{i,k=1}^n
\: E^{\alpha\beta} \cdot E^{\beta\gamma} @>\cong>> E^{\alpha\gamma}
$$
of vector bundle isomorphisms
$$
\phi^{\alpha\beta\gamma}_{ik} \: \bigoplus_{j=1}^n E^{\alpha\beta}_{ij}
\otimes E^{\beta\gamma}_{jk} @>\cong>> E^{\alpha\gamma}_{ik}
$$
over~$U_{\alpha\beta\gamma}$, for each triple $\alpha < \beta < \gamma$
in $\I$, such that
\item
the diagram
$$
\xymatrix{
E^{\alpha\beta} \cdot (E^{\beta\gamma} \cdot E^{\gamma\delta})
\ar[rr]^{\underline\alpha} \ar[d]_{id \cdot\phi^{\beta\gamma\delta}} &&
(E^{\alpha\beta} \cdot E^{\beta\gamma}) \cdot E^{\gamma\delta}
\ar[d]^{\phi^{\alpha\beta\gamma}\cdot id} \\
E^{\alpha\beta} \cdot E^{\beta\delta} \ar[r]_{\phi^{\alpha\beta\delta}} &
E^{\alpha\delta} &
E^{\alpha\gamma} \cdot E^{\gamma\delta} \ar[l]^{\phi^{\alpha\gamma\delta}}
}
$$
of vector bundle isomorphisms over~$U_{\alpha\beta\gamma\delta}$
commutes, for each chain $\alpha < \beta < \gamma < \delta$ in $\I$.
\endroster
Here $\underline\alpha$ denotes the (coherent) natural associativity
isomorphism for the matrix product $\cdot$ derived from the tensor
product $\otimes$ of vector bundles.  We call the $n \times n$ matrices
$E^{\alpha\beta}$ and $\phi^{\alpha\beta\gamma}$ the {\it gluing
bundles\/} and the {\it coherence isomorphisms\/} of the charted 2-vector
bundle $\E \downarrow X$, respectively.
\enddefinition

\definition{Definition 2.3}
Let $\E$ and $\F$ be two charted 2-vector bundles of rank~$n$ over~$X$, with
respect to the same ordered open cover $(\U, \I)$, with
gluing bundles $E^{\alpha\beta}$ and $F^{\alpha\beta}$ and coherence
isomorphisms $\phi^{\alpha\beta\gamma}$ and $\psi^{\alpha\beta\gamma}$,
respectively.  An {\it elementary change of trivializations}
$(T^\alpha, \tau^{\alpha\beta})$ from $\E$ to $\F$ is given by
\roster
\item
an $n \times n$ matrix $T^\alpha = (T^\alpha_{ij})_{i,j=1}^n$ of complex
vector bundles over~$U_\alpha$, for each $\alpha$ in $\I$, such that
over each point $x \in U_\alpha$ the integer matrix of fiber dimensions
$\dim(T^\alpha_x)$ has determinant $\pm1$, and
\item
an $n \times n$ matrix of vector bundle isomorphisms
$$
\tau^{\alpha\beta} = (\tau^{\alpha\beta}_{ij})_{i,j=1}^n
\: F^{\alpha\beta} \cdot T^\beta @>\cong>> T^\alpha \cdot E^{\alpha\beta}
$$
over~$U_{\alpha\beta}$, for each pair $\alpha < \beta$ in $\I$, such that
\item
the diagram
$$
\xymatrix{
F^{\alpha\beta} \cdot F^{\beta\gamma} \cdot T^\gamma
\ar[r]^{id \cdot \tau^{\beta\gamma}} \ar[d]_{\psi^{\alpha\beta\gamma} \cdot id} &
F^{\alpha\beta} \cdot T^\beta \cdot E^{\beta\gamma}
\ar[r]^{\tau^{\alpha\beta} \cdot id} &
T^\alpha \cdot E^{\alpha\beta} \cdot E^{\beta\gamma}
\ar[d]^{id \cdot \phi^{\alpha\beta\gamma}} \\
F^{\alpha\gamma} \cdot T^\gamma
\ar[rr]_{\tau^{\alpha\gamma}} &&
T^\alpha \cdot E^{\alpha\gamma}
}
$$
(natural associativity isomorphisms suppressed)
of vector bundle isomorphisms over~$U_{\alpha\beta\gamma}$ commutes,
for each triple $\alpha < \beta < \gamma$ in $\I$.
\endroster
\enddefinition

\definition{Definition 2.4}
Let $(\U, \I)$ and $(\U', \I')$ be two ordered open covers of $X$.
Suppose that there is an order-preserving {\it carrier function} $c \:
\I' \to \I$ such that for each $\alpha \in \I'$ there is an inclusion
$U'_\alpha \subset U_{c(\alpha)}$.  Then $(\U', \I')$ is a {\it
refinement\/} of $(\U, \I)$.

Let $\E$ be a charted 2-vector bundle of rank~$n$ over~$X$ with respect
to $(\U, \I)$, with gluing bundles $E^{\alpha\beta}$ and coherence
isomorphisms $\phi^{\alpha\beta\gamma}$.  Let
$$
c^*E^{\alpha\beta} = E^{c(\alpha)c(\beta)} | U'_{\alpha\beta}
$$
for $\alpha < \beta$ in $\I'$ and
$$
c^*\phi^{\alpha\beta\gamma} = \phi^{c(\alpha)c(\beta)c(\gamma)} |
U'_{\alpha\beta\gamma}
$$
for $\alpha < \beta < \gamma$ in $\I'$, be $n \times n$ matrices of
vector bundles and vector bundle isomorphisms over~$U'_{\alpha\beta}$
and $U'_{\alpha\beta\gamma}$, respectively.  Then there is a charted
2-vector bundle $c^*\E$ of rank~$n$ over~$X$ with respect to $(\U',
\I')$, with gluing bundles $c^*E^{\alpha\beta}$ and coherence isomorphisms
$c^*\phi^{\alpha\beta\gamma}$.  We say that $c^*\E$ is an {\it elementary
refinement\/} of $\E$.

More generally, two charted 2-vector bundles of rank~$n$ over~$X$ are
said to be {\it equivalent 2-vector bundles\/} if they can be linked by
a finite chain of elementary changes of trivializations and elementary
refinements.  (This is the notion of equivalence that appears to be
appropriate for our representability theorem~4.5.)
\enddefinition

\remark{Remark 2.5}
A charted 2-vector bundle of rank~$1$ consists of precisely the data
defining a {\it gerbe\/} over~$X$ with band~$\C^*$, as considered e.g.~by
Giraud \cite{Gi76}, Brylinski \cite{Br93} and Hitchin \cite{H99, \S1}.
There is a unitary form of the definition above, with Hermitian gluing
bundles and unitary coherence isomorphisms, and a unitary 2-vector bundle
of rank~$1$ is nothing but a gerbe with band $U(1)$.  In either case,
the set of equivalence classes of $\C^*$-gerbes or $U(1)$-gerbes over~$X$
is in natural bijection with the third integral cohomology group $H^3(X;
\Z)$ \cite{Br93, 5.2.10}.
\endremark

\definition{Definition 2.6}
Let $\E \downarrow X$ be a charted 2-vector bundle of rank~$n$, with
notation as above, and let $a \: Y \to X$ be a map of topological spaces.
Then there is a charted 2-vector bundle $a^*\E \downarrow Y$ of rank~$n$
obtained from $\E$ by {\it pullback\/} along $a$.  It is charted with
respect to the ordered open cover $(\U', \I)$ with $\U' = \{U'_\alpha
= a^{-1}(U_\alpha) \mid \alpha \in \I\}$.  It has gluing bundles
$a^*E^{\alpha\beta}$ obtained by pullback of the matrix of vector bundles
$E^{\alpha\beta}$ along $a \: U'_{\alpha\beta} \to U_{\alpha\beta}$, and
coherence isomorphisms $a^*\phi^{\alpha\beta\gamma}$ obtained by pullback
of the matrix of vector bundle isomorphisms $\phi^{\alpha\beta\gamma}$
along $a \: U'_{\alpha\beta\gamma} \to U_{\alpha\beta\gamma}$.
By definition there is then a {\it map\/} of charted 2-vector bundles
$\hat a \: a^*\E \to \E$ covering $a \: Y \to X$.
\enddefinition

\definition{Definition 2.7}
Let $\E \downarrow X$ and $\F \downarrow X$ be charted 2-vector
bundles with respect to the same ordered open cover $(\U, \I)$ of
$X$, of ranks~$n$ and~$m$, with gluing bundles $E^{\alpha\beta}$ and
$F^{\alpha\beta}$ and coherence isomorphisms $\phi^{\alpha\beta\gamma}$
and $\psi^{\alpha\beta\gamma}$, respectively.  Their {\it Whitney sum}
$\E \oplus \F \downarrow X$ is then the charted 2-vector bundle of rank
$(n+m)$ with gluing bundles given by the $(n+m) \times (n+m)$ matrix of
vector bundles
$$
\pmatrix E^{\alpha\beta} & 0 \\ 0 & F^{\alpha\beta} \endpmatrix
$$
and coherence isomorphisms given by the $(n+m) \times (n+m)$ matrix of
vector bundle isomorphisms
$$
\pmatrix \phi^{\alpha\beta\gamma} & 0 \\ 0 & \psi^{\alpha\beta\gamma}
\endpmatrix
\:
\pmatrix E^{\alpha\beta} & 0 \\ 0 & F^{\alpha\beta} \endpmatrix
\cdot
\pmatrix E^{\beta\gamma} & 0 \\ 0 & F^{\beta\gamma} \endpmatrix
@>\cong>>
\pmatrix E^{\alpha\gamma} & 0 \\ 0 & F^{\alpha\gamma} \endpmatrix
\,.
$$
There is an elementary change of trivializations from $\E \oplus \F$
to $\F \oplus \E$ given by the $(n+m) \times (n+m)$ matrix
$$
T^\alpha = \pmatrix 0 & I_m \\ I_n & 0 \endpmatrix
$$
for each $\alpha$ in $\I$, and identity isomorphisms $\tau^{\alpha\beta}$.
Here $I_n$ denotes the identity $n \times n$ matrix, with the trivial
rank~$1$ vector bundle in each diagonal entry and zero bundles elsewhere.
\enddefinition

\head 3. Algebraic $K$-theory of two-vector spaces \endhead

Let $(\B, \oplus, \otimes, \underline0, \underline1)$ be a {\it symmetric
bimonoidal category}, with sum and tensor functors
$$
\oplus, \otimes \: \B \times \B @>>> \B \,,
$$
and zero and unit objects $\underline0, \underline 1$ in $\B$.  These
satisfy associative, commutative and distributive laws, etc., up to a
list of natural isomorphisms, and these isomorphisms are {\it coherent\/}
in the sense that they fulfill a (long) list of compatibility conditions,
as presented by Laplaza in \cite{L72, \S1}.  We say that $\B$ is a {\it
bipermutative category\/} if the natural isomorphisms are almost all
identity morphisms, except for the commutative laws for $\oplus$ and
$\otimes$ and the left distributive law, and these in turn fulfill the
(shorter) list of compatibility conditions listed by May in \cite{Ma77,
\S VI.3}.

Suppose that $\B$ is {\it small}, i.e., that the class of objects of $\B$
is in fact a set.  Let $\pi_0(\B)$ be the set of path components of the
geometric realization of $\B$.  (Two objects of $\B$ are in the same
path component if and only if they can be linked by a finite chain of
morphisms in~$\B$.)  Then the sum and tensor functors induce sum and
product pairings that make $\pi_0(\B)$ into a commutative semi-ring
with zero and unit.  We can therefore think of the symmetric bimonoidal
category $\B$ as a kind of generalized commutative semi-ring.  Conversely,
any commutative semi-ring may be viewed as a discrete category, with
only identity morphisms, which is then a symmetric bimonoidal category.

The additive Grothendieck group completion $\Gr(\pi_0(\B))$ of the
commutative semi-ring $\pi_0(\B)$ is a commutative ring.  Likewise,
the geometric realization $|\B|$ can be group completed with respect
to the symmetric monoidal pairing induced by the sum functor $\oplus$,
and this group completion can take place at the categorical level, say
by Quillen's construction $\B^{-1}\B$ \cite{Gr76} or its generalization
$\B^+ = E\B \times_{\B} \B^2$ due to Thomason \cite{T79, 4.3.1}.  However,
the tensor functor $\otimes$ does not readily extend to $\B^{-1}\B$, as
was pointed out by Thomason \cite{T80}.  So $\B^{-1}\B$ is a symmetric
monoidal category, but usually not a symmetric bimonoidal category.

\example{Example 3.1}
Let $\V$ be the topological bipermutative category of finite dimensional
complex vector spaces, with set of objects $\N = \{0, 1, 2, \dots\}$
with $d \in \N$ interpreted as the complex vector space $\C^d$, and
morphism spaces
$$
\V(d, e) = \cases
U(d) & \text{if $d=e$,} \\ \emptyset & \text{otherwise}
\endcases
$$
from $d$ to $e$.  The sum functor $\oplus$ takes $(d, e)$ to $d+e$ and
embeds $U(d) \times U(e)$ into $U(d+e)$ by the block sum of matrices.
The tensor functor $\otimes$ takes $(d, e)$ to $de$ and maps $U(d)
\times U(e)$ to $U(de)$ by means of the left lexicographic ordering,
which identifies $\{1, \dots, d\} \times \{1, \dots, e\}$ with $\{1,
\dots, de\}$.  Both of these functors are continuous.
The zero and unit objects are $0$ and $1$, respectively.

In this case, the semi-ring $\pi_0(\V) = \N$ is that of the non-negative
integers, with additive group completion $\Gr(\N) = \Z$.  The geometric
realization $|\V| = \coprod_{d\ge0} BU(d)$ is the classifying space for
complex vector bundles, while its group completion $|\V^{-1}\V| \simeq
\Z \times BU$ classifies virtual vector bundles.  The latter space is
the infinite loop space underlying the spectrum $ku = \Spt(\V)$ that
represents connective complex topological $K$-theory, which is associated
to either of the symmetric monoidal categories $\V$ or $\V^{-1}\V$
by the procedure of Segal \cite{Se74}, as generalized by Shimada and
Shimakawa \cite{SS79} and Thomason \cite{T79, 4.2.1}.
\endexample

\definition{Definition 3.2}
Let $(\B, \oplus, \otimes, \underline0, \underline1)$ be a symmetric
bimonoidal category.  The category $M_n(\B)$ of {\it $n \times n$
matrices over~$\B$} has objects the matrices $V = (V_{ij})_{i,j=1}^n$
with entries that are objects of $\B$, and morphisms the matrices
$\phi = (\phi_{ij})_{i,j=1}^n$ with entries that are morphisms in $\B$.
The source (domain) of $\phi$ is the matrix of sources of the entries
$\phi_{ij}$, and similarly for targets (codomains).

There is a {\it matrix multiplication functor}
$$
M_n(\B) \times M_n(\B) @>\cdot>> M_n(\B)
$$
that takes two matrices $U = (U_{ij})_{i,j=1}^n$ and $V =
(V_{jk})_{j,k=1}^n$ to the matrix $W = U \cdot V = (W_{ik})_{i,k=1}^n$
with
$$
W_{ik} = \bigoplus_{j=1}^n U_{ij} \otimes V_{jk}
$$
for $i,k=1, \dots, n$.  In general, we need to make a definite choice
of how the $n$-fold sum is to be evaluated, say by bracketing from
the left.  When the direct sum functor is strictly associative, as in
the bipermutative case, the choice does not matter.

The {\it unit object} $I_n$ of $M_n(\B)$ is the $n \times n$
matrix with unit entries $\underline1$ on the diagonal and zero entries
$\underline 0$ everywhere else.
\enddefinition

\proclaim{Proposition 3.3}
$(M_n(\B), \cdot, I_n)$ is a monoidal category.
\endproclaim

In other words, the functor $\cdot$ is associative up to a natural
associativity isomorphism
$$
\underline\alpha \: U \cdot (V \cdot W) @>\cong>> (U \cdot V) \cdot W
$$
and unital with respect to $I_n$ up to natural left and right unitality
isomorphisms.  These are coherent, in the sense that they fulfill a
list of compatibility conditions, including the Mac\,Lane--Stasheff
pentagon axiom.  The proof of the proposition is a direct application
of Laplaza's first coherence theorem from \cite{L72, \S7}.

\definition{Definition 3.4}
Let $B$ be a commutative semi-ring with additive Grothendieck group
completion the commutative ring $A = \Gr(B)$.  Let $M_n(A)$ and $M_n(B)$
be the multiplicative monoids of $n \times n$ matrices with entries
in $A$ and $B$, respectively, and let $GL_n(A) \subset M_n(A)$ be the
subgroup of invertible $n \times n$ matrices with entries in $A$, i.e.,
those whose determinant is a unit in $A$.  Let the submonoid $GL_n(B)
\subset M_n(B)$ be the pullback in the diagram
$$
\xymatrix{
GL_n(B) \ar[r] \ar@{ >->}[d] &
GL_n(A) \ar@{ >->}[d] \\
M_n(B) \ar[r] &
M_n(A) \,.
}
$$
\enddefinition

\example{Example 3.5}
When $B = \N$ and $A = \Z$, $GL_n(\N) = M_n(\N) \cap GL_n(\Z)$ is the
monoid of $n \times n$ matrices with non-negative integer entries that
are invertible as integer matrices, i.e., have determinant $\pm1$.
It contains the elementary matrices that have entries $1$ on the
diagonal and in one other place, and $0$ entries elsewhere.
This is a larger monoid than the subgroup of units in $M_n(\N)$, which
only consists of the permutation matrices.
\endexample

\definition{Definition 3.6}
Let $\B$ be a symmetric bimonoidal category.  Let $GL_n(\B) \subset
M_n(\B)$ be the full subcategory with objects the matrices $V =
(V_{ij})_{i,j=1}^n$ whose matrix of path components $[V] =
([V_{ij}])_{i,j=1}^n$ lies in the submonoid $GL_n(\pi_0(\B)) \subset
M_n(\pi_0(\B))$.  We call $GL_n(\B)$ the category of {\it weakly
invertible $n \times n$ matrices over~$\B$}.
\enddefinition

\proclaim{Corollary 3.7}
$(GL_n(\B), \cdot, I_n)$ is a monoidal category.
\endproclaim

\definition{Definition 3.8}
Let $(\M, \cdot, e)$ be a monoidal category, and write $[p] = \{0 < 1 <
\dots < p\}$.  The {\it bar construction} $B\M$ is a simplicial category
$[p] \mapsto B_p\M$.  In simplicial degree $p$ the category $B_p\M$
has objects consisting of
\roster
\item
triangular arrays of objects $M^{\alpha\beta}$ of $\M$, for all $\alpha <
\beta$ in $[p]$, and
\item
isomorphisms
$$
\mu^{\alpha\beta\gamma} \: M^{\alpha\beta} \cdot M^{\beta\gamma}
@>\cong>> M^{\alpha\gamma}
$$
in $\M$, for all $\alpha < \beta < \gamma$ in $[p]$, such that
\item
the diagram of isomorphisms
$$
\xymatrix{
M^{\alpha\beta} \cdot (M^{\beta\gamma} \cdot M^{\gamma\delta})
\ar[rr]^{\underline\alpha} \ar[d]_{id \cdot \mu^{\beta\gamma\delta}} &&
(M^{\alpha\beta} \cdot M^{\beta\gamma}) \cdot M^{\gamma\delta}
\ar[d]^{\mu^{\alpha\beta\gamma} \cdot id} \\
M^{\alpha\beta} \cdot M^{\beta\delta} \ar[r]_{\mu^{\alpha\beta\delta}} &
M^{\alpha\delta} &
M^{\alpha\gamma} \cdot M^{\gamma\delta} \ar[l]^{\mu^{\alpha\gamma\delta}}
}
$$
commutes, for all $\alpha < \beta < \gamma < \delta$ in $[p]$.
\endroster
Here $\underline\alpha$ is the associativity isomorphism for the
monoidal operation $\cdot$ in $\M$.

The morphisms in $B_p\M$ from one object $(M_0^{\alpha\beta},
\mu_0^{\alpha\beta\gamma})$ to another $(M_1^{\alpha\beta},
\mu_1^{\alpha\beta\gamma})$ consist of a triangular array of morphisms
$\phi^{\alpha\beta} \: M_0^{\alpha\beta} \to M_1^{\alpha\beta}$ in
$\M$ for all $\alpha < \beta$ in $[p]$, such that the diagram
$$
\xymatrix{
M_0^{\alpha\beta} \cdot M_0^{\beta\gamma}
\ar[r]^-{\mu_0^{\alpha\beta\gamma}}
\ar[d]_{\phi^{\alpha\beta} \cdot \phi^{\beta\gamma}} &
M_0^{\alpha\gamma} \ar[d]^{\phi^{\alpha\gamma}} \\
M_1^{\alpha\beta} \cdot M_1^{\beta\gamma}
\ar[r]_-{\mu_1^{\alpha\beta\gamma}} &
M_1^{\alpha\gamma}
}
$$
commutes, for all $\alpha < \beta < \gamma$ in $[p]$.

To allow for degeneracy operators $f$ in the following paragraph,
let $M^{\alpha\alpha} = e$ be the unit object of $\M$, let
$\mu^{\alpha\alpha\beta}$ and $\mu^{\alpha\beta\beta}$ be the left
and right unitality isomorphisms for $\cdot$, respectively, and let
$\phi^{\alpha\alpha}$ be the identity morphism on $e$.

The simplicial structure on $B\M$ is given as follows.  For each
order-preserving function $f \: [q] \to [p]$ let the functor $f^* \: B_p\M
\to B_q\M$ take the object $(M^{\alpha\beta}, \mu^{\alpha\beta\gamma})$
of $B_p\M$ to the object of $B_q\M$ that consists of the triangular
array of objects $M^{f(\alpha)f(\beta)}$ for $\alpha < \beta$ in $[q]$
and the isomorphisms $\mu^{f(\alpha)f(\beta)f(\gamma)}$ for $\alpha <
\beta < \gamma$ in $[q]$.
\enddefinition

Each monoidal category $\M$ can be rigidified to an equivalent {\it
strict monoidal category} $\M_s$, i.e., one for which the associativity
isomorphism and the left and right unitality isomorphisms are all identity
morphisms \cite{ML71, XI.3.1}.  The usual strict bar construction for
$\M_s$ is a simplicial category $[p] \mapsto \M_s^p$, and corresponds
in simplicial degree~$p$ to the full subcategory of $B_p\M_s$ where all
the isomorphisms $\mu^{\alpha\beta\gamma}$ are identity morphisms.

\proclaim{Proposition 3.9}
The bar construction $B\M$ is equivalent to the strict bar construction
$[p] \mapsto \M_s^p$ for any strictly monoidal rigidification $\M_s$
of $\M$.
\endproclaim

This justifies calling $B\M$ the bar construction.  The proof is an
application of Quillen's theorem~A and the coherence theory for monoidal
categories.

\definition{Definition 3.10}
Let $\Ar\M = \Fun([1], \M)$ be the {\it arrow category\/} of $\M$, with
the morphisms of $\M$ as objects and commutative square diagrams in
$\M$ as morphisms.  There are obvious {\it source\/} and {\it target
functors} $s, t \: \Ar \M \to \M$.  Let $\Iso\M \subset \Ar\M$ be the
full subcategory with objects the isomorphisms of $\M$.
\enddefinition

\proclaim{Lemma 3.11}
Let $(\M, \cdot, e)$ be a monoidal category.
The category $B_2\M$ is the limit of the diagram
$$
\M \times \M @>\cdot>> \M @<s<< \Iso\M @>t>> \M \,.
$$
For $p\ge2$ each object or morphism of $B_p\M$ is uniquely determined
by the collection of its $2$-faces in $B_2\M$, which is indexed by the
set of monomorphisms $f \: [2] \to [p]$.
\endproclaim

Consider the symmetric bimonoidal category $\B$ as a kind of generalized
semi-ring.  The sum and tensor operations in $\B$ make the product
category $\B^n$ a generalized (right) module over~$\B$, for each
non-negative integer $n$.  The collection of $\B$-module homomorphisms
$\B^n \to \B^n$ is encoded in terms of (left) matrix multiplication by
the monoidal category $M_n(\B)$, and we shall interpret the monoidal
subcategory $GL_n(\B)$ as a category of {\it weak equivalences}
$\B^n @>\sim>> \B^n$.  This motivates the following definition.

\definition{Definition 3.12}
Let $\B$ be a symmetric bimonoidal category.  The {\it algebraic
$K$-theory\/} of the 2-category of (finitely generated free) modules
over~$\B$ is the loop space
$$
K(\B) = \Omega B \bigl( \coprod_{n\ge0} |BGL_n(\B)| \bigr) \,.
$$
Here $|BGL_n(\B)|$ is the geometric realization of the bar construction
on the monoidal category $GL_n(\B)$ of weakly invertible $n \times n$
matrices over~$\B$.  The block sum of matrices $GL_n(\B) \times GL_m(\B)
\to GL_{n+m}(\B)$ makes the coproduct $\coprod_{n\ge0} |BGL_n(\B)|$
a topological monoid.  The looped bar construction $\Omega B$ provides
a group completion of this topological monoid.

When $\B = \V$ is the category of finite dimensional complex vector
spaces, the (finitely generated free) modules over~$\V$ are called
{\it 2-vector spaces}, and $K(\V)$ is the algebraic $K$-theory of the
2-category of 2-vector spaces.

Let $GL_\infty(\B) = \colim_n GL_n(\B)$ be the infinite stabilization
with respect to block sum with the unit object in $GL_1(\B)$,
and write $B = \pi_0(\B)$ and $A = \Gr(B)$.  Then $K(\B) \simeq \Z
\times |BGL_\infty(\B)|^+$ by the McDuff--Segal group completion
theorem \cite{MS76}.  Here the superscript `$+$' refers to
Quillen's plus-construction with respect to the (maximal perfect)
commutator subgroup of $GL_\infty(A) \cong \pi_1 |BGL_\infty(\B)|$;
cf.~proposition~5.3 below.
\enddefinition

\head 4. Represented two-vector bundles \endhead

Let $X$ be a topological space, with an ordered open cover $(\U, \I)$.
Recall that all morphisms in $\V$ are isomorphisms, so $\Ar GL_n(\V) =
\Iso GL_n(\V)$.

\definition{Definition 4.1}
A {\it represented 2-vector bundle} $\E$ of rank~$n$ over~$X$ consists of
\roster
\item
a {\it gluing map}
$$
g^{\alpha\beta} \: U_{\alpha\beta} @>>> |GL_n(\V)|
$$
for each pair $\alpha < \beta$ in $\I$, and
\item
a {\it coherence map}
$$
h^{\alpha\beta\gamma} \: U_{\alpha\beta\gamma} @>>> |\!\Ar GL_n(\V)|
$$
satisfying $s \circ h^{\alpha\beta\gamma} = g^{\alpha\beta}
\cdot g^{\beta\gamma}$ and $t \circ h^{\alpha\beta\gamma} =
g^{\alpha\gamma}$ over~$U_{\alpha\beta\gamma}$,
for each triple $\alpha < \beta < \gamma$ in $\I$, 
such that
\item
the {\it $2$-cocycle condition}
$$
h^{\alpha\gamma\delta} \circ (h^{\alpha\beta\gamma} \cdot id)
\circ \underline\alpha =
h^{\alpha\beta\delta} \circ (id \cdot h^{\beta\gamma\delta})
$$
holds over~$U_{\alpha\beta\gamma\delta}$ for all
$\alpha < \beta < \gamma < \delta$ in $\I$.
\endroster
\enddefinition

There is a suitably defined notion of {\it equivalence\/} of represented
2-vector bundles, which we omit to formulate here, but cf.~definitions~2.3
and~2.4.

\definition{Definition 4.2}
Let $E(d) = EU(d) \times_{U(d)} \C^d \downarrow BU(d)$ be the universal
$\C^d$-bundle over~$BU(d)$.  There is a universal $n \times n$ matrix
$$
E = (E_{ij})_{i,j=1}^n
$$
of Hermitian vector bundles over~$|GL_n(\V)|$.  Over the path component
$|GL_n(\V)_D| = \prod_{i,j=1}^n BU(d_{ij})$ for $D = (d_{ij})_{i,j=1}^n$
in $GL_n(\N)$, the $(i,j)$-th entry in $E$ is the pullback of the
universal bundle $E(d_{ij})$ along the projection $|GL_n(\V)_D| \to
BU(d_{ij})$.

Let $|\!\Ar U(d)|$ be the geometric realization of the arrow category $\Ar
U(d)$, where $U(d)$ is viewed as a topological groupoid with one object.
Each pair $(A, B) \in U(d)^2$ defines a morphism from $C$ to $(A, B) \cdot
C = B C A^{-1}$, so
$$
|\!\Ar U(d)| \cong EU(d)^2 \times_{U(d)^2} U(d)
$$
equals the Borel construction for this (left) action of $U(d)^2$
on $U(d)$.  There are source and target maps $s, t \: |\!\Ar U(d)| \to
BU(d)$, which take the $1$-simplex represented by a morphism $(A, B)$ to
the $1$-simplices represented by the morphisms $A$ and $B$, respectively.
By considering each element in $U(d)$ as a unitary isomorphism $\C^d \to
\C^d$ one obtains a universal unitary vector bundle isomorphism $\phi(d)
\: s^*E(d) @>\cong>> t^*E(d)$

There is a universal $n \times n$ matrix of unitary vector bundle
isomorphisms
$$
\phi \: s^*E \cong t^*E
$$
over~$|\!\Ar GL_n(\V)|$.  Over the path component $|\!\Ar GL_n(\V)_D| =
\prod_{i,j=1}^n |\!\Ar U(d_{ij})|$ for $D$ as above, the $(i,j)$-th entry
in $\phi$ is the pullback of the universal isomorphism $\phi(d_{ij})$
along the projection $|\!\Ar GL_n(\V)_D| \to |\!\Ar U(d_{ij})|$.
\enddefinition

\proclaim{Lemma 4.3}
Let $\E$ be a represented 2-vector bundle with gluing maps $g^{\alpha\beta}$
and coherence maps $h^{\alpha\beta\gamma}$.  There is an associated
charted 2-vector bundle with gluing bundles
$$
E^{\alpha\beta} = (g^{\alpha\beta})^*(E)
$$
over~$U_{\alpha\beta}$ and coherence isomorphisms
$$
\phi^{\alpha\beta\gamma} = (h^{\alpha\beta\gamma})^*(\phi)
\: E^{\alpha\beta} \cdot E^{\beta\gamma} =
(g^{\alpha\beta} \cdot g^{\beta\gamma})^*(E) @>\cong>>
(g^{\alpha\gamma})^*(E) = E^{\alpha\gamma}
$$
over~$U_{\alpha\beta\gamma}$.  The association induces a bijection
between the equivalence classes of represented 2-vector bundles and the
equivalence classes of charted 2-vector bundles of rank~$n$ over~$X$.
\endproclaim

\definition{Definition 4.4}
Let $\twoVect_n(X)$ be the set of equivalence classes of 2-vector
bundles of rank~$n$ over~$X$.  For path-connected $X$ let $\twoVect(X)
= \coprod_{n\ge0} \twoVect_n(X)$.  Whitney sum (definition~2.7) defines
a pairing that makes $\twoVect(X)$ an abelian monoid.
\enddefinition

\proclaim{Theorem 4.5}
Let $X$ be a finite CW complex.  There are natural bijections
$$
\twoVect_n(X) \cong [X, |BGL_n(\V)|]
$$
and
$$
\twoVect(X) \cong [X, \coprod_{n\ge0} |BGL_n(\V)|] \,.
$$
\endproclaim

To explain the first correspondence, from which the second follows,
we use the following construction.

\definition{Definition 4.6}
Let $(\U, \I)$ be an ordered open cover of $X$.  The {\it Mayer--Vietoris
blow-up} $MV(\U)$ of $X$ with respect to $\U$ is the simplicial space
with $p$-simplices
$$
MV_p(\U) = \coprod_{\alpha_0 \le \dots \le \alpha_p} U_{\alpha_0\dots\alpha_p}
$$
with $\alpha_0 \le \dots \le \alpha_p$ in $\I$.  The $i$-th simplicial
face map is a coproduct of inclusions $U_{\alpha_0\dots\alpha_p}
\subset U_{\alpha_0\dots\hat\alpha_i\dots\alpha_p}$, and similarly for
the degeneracy maps.  The inclusions $U_{\alpha_0\dots\alpha_p} \subset
X$ combine to a natural map $e \: |MV(\U)| \to X$, which is a (weak)
homotopy equivalence.
\enddefinition

\demo{Sketch of proof of theorem~4.5}
By lemma~3.11, a simplicial map $g \: MV(\U) \to |BGL_n(\V)|$ is uniquely
determined by its components in simplicial degrees $1$ and $2$.  The first
of these is a map
$$
g_1 \: MV_1(\U) = \coprod_{\alpha \le
\beta} U_{\alpha\beta} @>>> |B_1GL_n(\V)| = |GL_n(\V)|
$$
which is a coproduct of gluing maps $g^{\alpha\beta} \: U_{\alpha\beta}
\to |GL_n(\V)|$.  The second is a map
$$
g_2 \: MV_2(\U) = \coprod_{\alpha \le \beta \le \gamma} U_{\alpha\beta\gamma}
\to |B_2GL_n(\V)| \,.
$$
The simplicial identities and lemma~3.11 imply that $g_2$ is determined
by $g_1$ and a coproduct of coherence maps $h^{\alpha\beta\gamma} \:
U_{\alpha\beta\gamma} \to |\!\Ar GL_n(\V)|$.  Hence such a simplicial map
$g$ corresponds bijectively to a represented 2-vector bundle of rank~$n$
over~$X$.

Any map $f \: X \to |BGL_n(\V)|$ can be composed with the weak equivalence
$e \: |MV(\U)| \to X$ to give a map of spaces $fe \: |MV(\U)| \to
|BGL_n(\V)|$, which is homotopic to a simplicial map $g$ if $\U$ is
a good cover, and for reasonable $X$ any open ordered cover can be refined to a
good one.  The homotopy class of $f$ corresponds to the equivalence class
of the represented 2-vector bundle determined by the simplicial map~$g$.
\qed
\enddemo

\remark{Remark 4.7}
We wish to interpret the 2-vector bundles over~$X$ as (effective)
$0$-cycles for some cohomology theory at $X$.  Such theories are
group-valued, so a first approximation to the $0$-th cohomology group
at $X$ could be the Grothendieck group $\Gr(\twoVect(X))$ of formal
differences of 2-vector bundles over~$X$.  The analogous construction
for ordinary vector bundles works well to define topological $K$-theory,
but for 2-vector bundles this algebraically group completed functor is
not even representable, like in the case of the algebraic $K$-theory of
a discrete ring.  We thank Haynes Miller for reminding us of this issue.

Instead we follow Quillen and perform the group completion at the space
level, which leads to the algebraic $K$-theory space
$$
\align
K(\V) &= \Omega B \bigl( \coprod_{n\ge0} |BGL_n(\V)| \bigr) \\
&\simeq \Z \times |BGL_\infty(\V)|^+
\endalign
$$
from definition~3.12.  But what theory does this loop space represent?
One interpretation is provided by the theory of virtual flat fibrations,
presented by Karoubi in \cite{K87, Ch.~III}, leading to what we shall
call virtual 2-vector bundles.  Another interpretation could be given
using the homology bordism theory of Hausmann and Vogel \cite{HV78}.
\endremark

\definition{Definition 4.8}
Let $X$ be a space.  An {\it acyclic fibration\/} over~$X$ is a Serre
fibration $a \: Y \to X$ such that the homotopy fiber at each point $x \in
X$ has the integral homology of a point, i.e., $\tilde H_*(\hofib_x(a);
\Z) = 0$.  A map of acyclic fibrations from $a' \: Y' \to X$ to $a \:
Y \to X$ is a map $f \: Y' \to Y$ with $af = a'$.

A {\it virtual 2-vector bundle\/} over~$X$ is described by an acyclic
fibration $a \: Y \to X$ and a 2-vector bundle $\E \downarrow Y$.
We write $\E \downarrow Y @>a>> X$.  Given a map $f \: Y' \to Y$ of
acyclic fibrations over~$X$ there is an induced 2-vector bundle $f^*\E
\downarrow Y'$.  The virtual 2-vector bundles described by $\E \downarrow
Y @>a>> X$ and $f^*\E \downarrow Y' @>a'>> X$ are declared to be {\it
equivalent\/} as virtual 2-vector bundles over~$X$.
\enddefinition

\proclaim{Lemma 4.9}
The abelian monoid of equivalence classes of virtual 2-vector bundles
over~$X$ is the colimit
$$
\colim_{a \: Y \to X} \twoVect(Y)
$$
where $a \: Y \to X$ ranges over the category of acyclic fibrations
over~$X$.  Its Grothendieck group completion is isomorphic to the colimit
$$
\colim_{a \: Y \to X} \Gr(\twoVect(Y)) \,.
$$
\endproclaim

The functor $Y \mapsto \twoVect(Y)$ factors through the homotopy category
of acyclic fibrations over $X$, which is directed.

The following result says that formal differences of virtual 2-vector
bundles over~$X$ are the geometric objects that constitute cycles for the
contravariant homotopy functor represented by the algebraic $K$-theory
space $K(\V)$.  Compare \cite{K87, III.3.11}.  So $K(\V)$ represents
sheaf cohomology for the topology of acyclic fibrations, with coefficients
in the abelian presheaf $Y \mapsto \Gr(\twoVect(Y))$ given by
the Grothendieck group completion of the abelian monoid of
equivalence classes of 2-vector bundles.

\proclaim{Theorem 4.10}
Let $X$ be a finite CW complex.
There is a natural group isomorphism
$$
\colim_{a \: Y \to X} \Gr(\twoVect(Y)) \cong [X, K(\V)]
$$
where $a \: Y \to X$ ranges over the category of acyclic fibrations over
$X$.  Restricted to $\Gr(\twoVect(X))$ (with $a = id$) the isomorphism
extends the canonical monoid homomorphism $\twoVect(X) \cong [X,
\coprod_{n\ge0} |BGL_n(\V)|] \to [X, K(\V)]$.
\endproclaim

\remark{Remark 4.11}
The passage to sheaf cohomology would be unnecessary if we replaced
$\V$ by a different symmetric bimonoidal category $\B$ such that
each $\pi_0(GL_n(\B))$ is abelian.  This might entail an extension of
the category of vector spaces to allow generalized vector spaces of
arbitrary real, or even complex, dimension, parallel to the inclusion
of the integers into the real or complex numbers.  Such an extension is
reminiscent of a category of representations of a suitable $C^*$-algebra,
but we know of no clear interpretation of this approach.
\endremark

\head 5. Algebraic $K$-theory of topological $K$-theory \endhead

Is the contravariant homotopy functor $X \mapsto [X, K(\V)] = K(\V)^0(X)$
part of a cohomology theory, and if so, what is the spectrum representing
that theory?

The topological symmetric bimonoidal category $\V$ plays the role
of a generalized commutative semi-ring in our definition of $K(\V)$.
Its additive group completion $\V^{-1}\V$ correspondingly plays the role
of a generalized commutative ring.  This may be tricky to realize at the
level of symmetric bimonoidal categories, but the connective topological
$K$-theory spectrum $ku = \Spt(\V)$ associated to the additive topological
symmetric monoidal structure of $\V$ is an $E_\infty$ ring spectrum,
and hence a commutative algebra over the sphere spectrum~$\SS$.

The algebraic $K$-theory of an $\SS$-algebra $\A$ can on one hand be
defined as the Waldhausen algebraic $K$-theory \cite{W85} of a category
with cofibrations and weak equivalences, with objects the finite cell
$\A$-modules, morphisms the $\A$-module maps and weak equivalences
the stable equivalences.  Alternatively, it can be defined as a group
completion
$$
K(\A) = \Omega B \bigl( \coprod_{n\ge0} B\widehat{GL}_n(\A) \bigr)
$$
where $\widehat{GL}_n(\A)$ is essentially the topological monoid of
$\A$-module maps $\A^n \to \A^n$ that are stable equivalences.  The former
definition produces a spectrum, so the space $K(\A)$ is in fact an
infinite loop space, and its deloopings represent a cohomology theory.

The passage from modules over the semi-ring object $\V$ to modules
over the ring object $ku$ corresponds to maps $|GL_n(\V)| \to
\widehat{GL}_n(ku)$ and a map $K(\V) \to K(ku)$.

\proclaim{Conjecture 5.1}
There is a weak equivalence $K(\V) \simeq K(ku)$.  More generally,
$K(\B) \simeq K(\A)$ for each symmetric bimonoidal category $\B$ with
associated commutative $\SS$-algebra $\A = \Spt(\B)$
\endproclaim

\remark{Remark 5.2}
The conjecture asserts that the contravariant homotopy functor $X
\mapsto [X, K(\V)]$ with $0$-cycles given by the virtual 2-vector
bundles over~$X$ is the $0$-th cohomology group for the cohomology
theory represented by the spectrum $K(ku)$ given by the algebraic
$K$-theory of connective topological $K$-theory.  We consider the
virtual 2-vector bundles over~$X$ to be sufficiently geometric objects
(like complex vector bundles), that this cohomology theory then admits
as geometric an interpretation as the classical examples of de\,Rham
cohomology, topological $K$-theory and complex bordism.
\endremark

\medskip

As a first (weak) justification of this conjecture, recall that to
the eyes of algebraic $K$-theory the block sum operation $(g, h)
\mapsto \left[\smallmatrix g & 0 \\ 0 & h \endsmallmatrix\right]$ is
identified with the stabilized matrix multiplication $(g, h) \mapsto
\left[\smallmatrix gh & 0 \\ 0 & I \endsmallmatrix\right]$, where $I$ is
an identity matrix.  The group completion in the definition of algebraic
$K$-theory adjoins inverses to the block sum operation, and thus also to
the stabilized matrix multiplication.  In particular, for each elementary
$n \times n$ matrix $e_{ij}(V)$ over~$\B$ with $(i,j)$-th off-diagonal
entry equal to an object $V$ of $\B$, the inverse matrix $e_{ij}(-V)$
is formally adjoined as far as algebraic $K$-theory is concerned.
Hence the formal negatives $-V$ in $\B^{-1}\B$ are already present,
in this weak sense.

A stronger indication that the conjecture should hold true is provided by
the following special case.  Recall that a commutative semi-ring is the
same as a (small) symmetric bimonoidal category that is {\it discrete},
i.e., has only identity morphisms.

\proclaim{Proposition 5.3}
Let $B$ be a commutative semi-ring, with additive Grothendieck
group completion $A = \Gr(B)$.  The semi-ring homomorphism $B \to A$
induces a weak equivalence
$$
BGL_\infty(B) @>\simeq>> BGL_\infty(A)
$$
and thus a weak equivalence $K(B) \simeq K(A)$.  In particular, there is
a weak equivalence $K(\N) \simeq K(\Z)$.
\endproclaim

A proof uses the following application of Quillen's theorem~B.

\proclaim{Lemma 5.4}
Let $f \: M \to G$ be a monoid homomorphism from a monoid $M$ to a group
$G$.  Write $mg = f(m) \cdot g$. Let $Q = B(*, M, G)$ be the category
with objects $g \in G$ and morphisms $(m, g) \in M \times G$ from $mg$ to
$g$:
$$
mg @>(m, g)>> g \,.
$$
Then there is a fiber sequence up to homotopy
$$
|Q| @>>> BM @>Bf>> BG \,.
$$
\endproclaim

\demo{Sketch of proof of proposition~5.3}
Applying lemma~5.4 to the monoids $M_n = GL_n(B)$ and groups $G_n =
GL_n(A)$ we obtain categories $Q_n$ for each natural number~$n$.
There are stabilization maps $i \: Q_n \to Q_{n+1}$, $M_n \to M_{n+1}$
and $G_n \to G_{n+1}$, with (homotopy) colimits $Q_\infty$, $M_\infty$
and $G_\infty$, and a quasi-fibration
$$
|Q_\infty| @>>> BGL_\infty(B) @>>> BGL_\infty(A) \,.
$$
It suffices to show that each stabilization map $i \: |Q_n| \to
|Q_{n+1}|$ is weakly null-homotopic, because then $|Q_\infty|$ is weakly
contractible.

For each full subcategory $K \subset Q_n$ with finitely many objects,
the restricted stabilization functor $i|K$ takes $g$ to $i(g) =
\left[\smallmatrix g & 0 \\ 0 & 1 \endsmallmatrix\right]$.  It receives a
natural transformation from a functor $j \: K \to Q_{n+1}$ that maps $g$
to $j(g) = \left[\smallmatrix g & v \\ 0 & 1 \endsmallmatrix\right]$ for
some column vector $v = v(g)$ with positive entries in $B$.  The trick is
to construct $v(g)$ inductively for the finite set of objects $g$ of $K$,
so that $v(mg)$ is sufficiently positive with respect to $m \cdot v(g)$
for all morphisms $mg \to g$ in $K$.

Furthermore, the finiteness of $K$ ensures that there is a row vector
$w$ with entries in $B$ and an object $h = \left[\smallmatrix
I_n & 0 \\ -w & 1 \endsmallmatrix\right]$ of $Q_{n+1}$ such that there
is a natural transformation from $j$ to the constant functor to $h$.
These two natural transformations provide a homotopy from $i|K$ to a
constant map.  As $K$ was arbitrary with finitely many objects, this
means that $i$ is weakly null-homotopic.
\qed
\enddemo

\remark{Remark 5.5}
If there exists a symmetric bimonoidal category $\W$ and a functor $\V
\to \W$ of symmetric bimonoidal categories that induces an additive
equivalence from $\V^{-1}\V$ to $\W$, then most likely the line of
argument sketched above in the case of commutative semi-rings can be
adapted to the symmetric bimonoidal situation.  This could provide one
line of proof toward conjecture~5.1.  Similar remarks apply for a general
symmetric bimonoidal category $\B$ in place of $\V$.
\endremark

\head 6. Forms of elliptic cohomology \endhead

In this section we shall view the algebraic $K$-theory $K(\A)$ of an
$\SS$-algebra $\A$ as a spectrum, rather than as a space.

We shall argue that the algebraic $K$-theory $K(ku)$ of the connective
topological $K$-theory spectrum $ku$ is a connective form of elliptic
cohomology, in the sense that it detects homotopy theoretic phenomena
related to $v_2$-periodicity, much like how topological $K$-theory
detects phenomena related to $v_1$-periodicity (which is really the same
as Bott periodicity) and how rational cohomology detects phenomena related
to $v_0$-periodicity.  Furthermore, from this point of view the homotopy
type of $K(ku)$ is robust with respect to changes in the interpretation
of the phrase ``algebraic $K$-theory of topological $K$-theory''.

We first introduce a filtration of the class of spectra that is related to
the {\it chromatic filtration\/} given by the property of being Bousfield
local with respect to some Johnson--Wilson theory $E(n)$ (cf.~Ravenel
\cite{Ra92, \S7}), but is more appropriate for the connective spectra
that arise from algebraic $K$-theory.  Our notion is also more closely
linked to aspects of $v_n$-periodicity than to being $E(n)$-local.

Let $p$ be a prime, $K(n)$ the $n$-th Morava $K$-theory at $p$ and $F$
a $p$-local finite CW spectrum.  The least number $0 \le n \le \infty$
such that $K(n)_*(F)$ is non-trivial is called the {\it chromatic
type\/} of $F$.  By the Hopkins--Smith periodicity theorem \cite{HS98,
Thm.~9}, $F$ admits a {\it $v_n$-self map} $v \: \Sigma^d F \to F$ such
that $K(m)_*(v)$ is an isomorphism for $m = n$ and zero for $m \ne n$.
The $v_n$-self map is sufficiently unique for the mapping telescope
$$
v^{-1}F = \Tel \bigl( F @>v>> \Sigma^{-d}F @>v>> \dots \bigr)
$$
to be well-defined up to homotopy.  The class of all $p$-local finite CW
spectra of chromatic type $\ge n$ is closed under weak equivalences
and the formation of homotopy cofibers, desuspensions and retracts, so
we say that the full subcategory that it generates is a {\it thick\/}
subcategory.  By the Hopkins--Smith thick subcategory theorem \cite{HS98,
Thm.~7}, any thick subcategory of the category of $p$-local finite CW
spectra has this precise form, for a unique number $0 \le n \le \infty$.

\definition{Definition 6.1}
Let $X$ be a spectrum, and let $\T_X$ be the full subcategory of $p$-local
finite CW spectra $F$ for which the localization map
$$
F \wedge X @>>> v^{-1}F \wedge X
$$
induces an isomorphism on homotopy groups in all sufficiently high
degrees.  Then $\T_X$ is a thick subcategory, hence consists of the
spectra $F$ of chromatic type $\ge n$ for some unique number $0 \le n
\le \infty$.  We call this number $n = \telecom(X)$ the {\it telescopic
complexity\/} of $X$.  (This abbreviation is due to Matthew Ando.)
\enddefinition

\proclaim{Lemma 6.2}
If $Y$ is the $k$-connected cover of $X$, for some integer $k$,
then $X$ and $Y$ have the same telescopic complexity.

Let $X \to Y \to Z$ be a cofiber sequence and $m = \max\{\telecom(X),
\telecom(Y)\}$.  If $\telecom(X) \ne \telecom(Y)$ then $\telecom(Z) =
m$, otherwise $\telecom(Z) \le m$.

If $Y$ is a (de-)suspension of $X$ then $X$ and $Y$ have the same
telescopic complexity.

If $Y$ is a retract of $X$ then $\telecom(Y) \le \telecom(X)$.

If $X$ is an $E(n)$-local spectrum then $X$ has telescopic complexity
$\le n$.
\endproclaim

\example{Examples 6.3}
(1) Integral, rational, real and complex cohomology ($H\Z$, $H\Q$, $H\R$
or $H\C$) all have telescopic complexity~$0$.

(2) Connective or periodic, real or complex topological $K$-theory ($ko$,
$ku$, $KO$ or $KU$) all have telescopic complexity~$1$.  The {\'e}tale
$K$-theory $K^{et}(R)$ of a ring $R = \Cal O_{F, S}$ of $S$-integers in a
local or global number field has telescopic complexity~$1$, and so does
the algebraic $K$-theory $K(R)$ if the Lichtenbaum--Quillen conjecture
holds for the ring $R$.

(3) An Ando--Hopkins--Strickland \cite{AHS01} {\it elliptic spectrum} $(E,
C, t)$ has telescopic complexity $\le 2$, and the telescopic complexity
equals~$2$ if and only if the elliptic curve $C$ over~$R = \pi_0(E)$
has a supersingular specialization over some point of $\Spec(R)$.

(4) The Hopkins--Mahowald--Miller {\it topological modular forms\/}
spectra $\tmf$ and $\TMF$ have telescopic complexity~$2$.

(5) The Johnson--Wilson spectrum $E(n)$ and its connective form, the
Brown--Peterson spectrum $BP\langle n\rangle$, both have telescopic
complexity~$n$.

(6) The sphere spectrum $\SS$ and the complex bordism spectrum $MU$
have infinite telescopic complexity.
\endexample

Let $V(1)$ be the four-cell Smith--Toda spectrum with $BP_*(V(1)) =
BP_*/(p, v_1)$.  For $p\ge5$ it exists as a commutative ring spectrum.
It has chromatic type~$2$, and there is a $v_2$-self map $v \:
\Sigma^{2p^2-2} V(1) \to V(1)$ inducing multiplication by the class $v_2
\in \pi_{2p^2-2} V(1)$.  We write $V(1)_*(X) = \pi_*(V(1) \wedge X)$
for the $V(1)$-homotopy groups of $X$, which are naturally a graded
module over $P(v_2) = \Bbb F_p[v_2]$.

Let $X_{(p)}$ and $X^\wedge_p$ denote the $p$-localization and
$p$-completion of a spectrum $X$, respectively.  The first Brown--Peterson
spectrum $\ell = BP\langle1\rangle$ is the connective $p$-local Adams
summand of $ku_{(p)}$, and its $p$-completion $\ell^\wedge_p$ is the
connective $p$-complete Adams summand of $ku^\wedge_p$.  These are all
known to be commutative $\SS$-algebras.

The spectrum $TC(\ell^\wedge_p)$ occurring in the following statement is
the {\it topological cyclic homology\/} of $\ell^\wedge_p$, as defined
by B{\"o}kstedt, Hsiang and Madsen \cite{BHM93}.  The theorem is proved
in \cite{AR02, 0.3} by an elaborate but explicit computation of its
$V(1)$-homotopy groups, starting from the corresponding $V(1)$-homotopy
groups of the topological Hochschild homology $THH(\ell^\wedge_p)$.

\proclaim{Theorem 6.4 (Ausoni--Rognes)}
Let $p\ge5$.  The algebraic $K$-theory spectrum $K(\ell^\wedge_p)$
of the connective $p$-complete Adams summand $\ell^\wedge_p$ has
telescopic complexity~$2$.  More precisely, there is an exact sequence
of $P(v_2)$-modules
$$
0 \to \Sigma^{2p-3} \Bbb F_p \to V(1)_* K(\ell^\wedge_p) @>\trc>>
V(1)_* TC(\ell^\wedge_p) \to \Sigma^{-1} \Bbb F_p \to 0
$$
and an isomorphism of $P(v_2)$-modules
$$
\multline
V(1)_* TC(\ell^\wedge_p) \cong P(v_2) \otimes \\
\bigl(\,
E(\partial, \lambda_1, \lambda_2) \oplus
E(\lambda_2) \{\lambda_1 t^d \mid 0 < d < p\} \oplus
E(\lambda_1) \{\lambda_2 t^{dp} \mid 0 < d < p\}
\,\bigr) \,.
\endmultline
$$
Here $\partial$, $t$, $\lambda_1$ and $\lambda_2$ have degrees
$-1$, $-2$, $2p-1$ and $2p^2-1$, respectively.  Hence $V(1)_*
TC(\ell^\wedge_p)$ is free of rank~$(4p+4)$ over~$P(v_2)$, and agrees
with its $v_2$-localization in sufficiently high degrees.
\endproclaim

Since $K(\ell^\wedge_p)$ has telescopic complexity~$2$, it has a chance to
detect $v_2$-periodic families in $\pi_* V(1)$.  This is indeed the case.
Let $\alpha_1 \in \pi_{2p-3} V(1)$ and $\beta_1' \in \pi_{2p^2-2p-1}
V(1)$ be the classes represented in the Adams spectral sequence by the
cobar $1$-cycles $h_{10} = [\bar\xi_1]$ and $h_{11} = [\bar\xi_1^p]$,
respectively.
There are maps $V(1) \to v_2^{-1} V(1) \to L_2 V(1)$, and Ravenel
\cite{Ra86, 6.3.22} computed
$$
\pi_* L_2 V(1) \cong P(v_2, v_2^{-1}) \otimes E(\zeta) \{1,
h_{10}, h_{11}, g_0, g_1, h_{11}g_0 = h_{10}g_1\}
$$
for $p\ge5$.  Hence $\pi_* L_2 V(1)$ contains twelve $v_2$-periodic
families.  The telescope conjecture asserted that $v_2^{-1} V(1) \to
L_2 V(1)$ might be an equivalence, but this is now considered to be
unlikely \cite{MRS01}.  The following detection result can be read off
from \cite{AR02, 4.8}, and shows that $K(\ell^\wedge_p)$ detects the same
kind of homotopy theoretic phenomena as $E(2)$ or an elliptic spectrum.

\proclaim{Proposition 6.5}
The unit map $\SS \to K(\ell^\wedge_p)$ induces a $P(v_2)$-module homomorphism
$\pi_* V(1) \to V(1)_* K(\ell^\wedge_p)$ which takes $1$, $\alpha_1$ and
$\beta_1'$ to $1$, $t\lambda_1$ and $t^p\lambda_2$, respectively.
Hence $V(1)_* K(\ell^\wedge_p)$ detects the $v_2$-periodic families in $\pi_*
V(1)$ generated by these three classes.
\endproclaim

Turning to the whole connective $p$-complete topological $K$-theory
spectrum $ku^\wedge_p$, there is a map $\ell^\wedge_p \to ku^\wedge_p$ of
commutative $\SS$-algebras.  It induces a natural map $K(\ell^\wedge_p)
\to K(ku^\wedge_p)$, and there is a transfer map $K(ku^\wedge_p) \to
K(\ell^\wedge_p)$ such that the composite self-map of $K(\ell^\wedge_p)$
is multiplication by $(p-1)$.  Hence the composite map is a $p$-local
equivalence.

\proclaim{Lemma 6.6}
The algebraic $K$-theory spectrum $K(ku^\wedge_p)$ of connective
$p$-complete topological $K$-theory $ku^\wedge_p$ contains
$K(\ell^\wedge_p)$ as a $p$-local retract, hence has telescopic complexity
$\ge 2$.
\endproclaim

Most likely $K(ku^\wedge_p)$ also has telescopic complexity exactly~$2$.
It may be possible to prove this directly by computing $V(1)_*
TC(ku^\wedge_p)$, by similar methods as in \cite{AR02}, but the algebra
involved for $ku^\wedge_p$ is much more intricate than it was for the
Adams summand.  Some progress in this direction has recently been made
by Ausoni.

The following consequence of a theorem of the second author \cite{D97,
p.~224} allows us to compare the algebraic $K$-theory of $ku^\wedge_p$
to that of the integral spectra $ku$ and $K(\C)$.

\proclaim{Theorem 6.7 (Dundas)}
Let $\A$ be a connective $\SS$-algebra.  The commutative square
$$
\xymatrix{
K(\A) \ar[r] \ar[d] &
K(\A^\wedge_p) \ar[d] \\
K(\pi_0(\A)) \ar[r] &
K(\pi_0(\A^\wedge_p))
}
$$
becomes homotopy Cartesian after $p$-completion.
\endproclaim

We apply this with $\A = ku$ or $\A = K(\C)$.  Also in the second case
$\A^\wedge_p \simeq ku^\wedge_p$, by Suslin's theorem on the algebraic
$K$-theory of algebraically closed fields \cite{Su83}.  Then $\pi_0(\A)
= \Z$ and $\pi_0(\A^\wedge_p) = \Z_p$.  It is known that $K(\Z_p)$
has telescopic complexity~$1$, by B{\"o}kstedt--Madsen \cite{BM94},
\cite{BM95} for $p$ odd and by the third author \cite{Ro99} for $p=2$.
It is also known that $K(\Z)$ has telescopic complexity~$1$ for $p=2$, by
Voevodsky's proof of the Milnor conjecture and Rognes--Weibel \cite{RW00}.
For $p$ odd it would follow from the Lichtenbaum--Quillen conjecture
for $K(\Z)$ at~$p$ that $K(\Z)$ has telescopic complexity~$1$, and this
now seems to be close to a theorem by the work of Voevodsky, Rost and
Positselski.

\proclaim{Proposition 6.8}
Suppose that $K(\Z)$ has telescopic complexity~$1$ at a prime $p\ge5$.
Then $K(ku)$ and $K(K(\C))$ have the same telescopic complexity as
$K(ku^\wedge_p)$, which is $\ge2$.
\endproclaim

More generally it is natural to expect that $K(K(R))$ has telescopic
complexity~$2$ for each ring of $S$-integers $R = \Cal O_{F, S}$ in a
local or global number field $F$, including the initial case $K(K(\Z))$.
A discussion of such a conjecture has been given in lectures by the
third author, but should take place in the context of {\'e}tale covers
or Galois extensions of commutative $\SS$-algebras, which would take us
too far afield here.

The difference between the connective and periodic topological $K$-theory
spectra $ku$ and $KU$ may also not affect their algebraic $K$-theories
greatly.  There is a (localization) fiber sequence
$$
K(\eusm C^{KU}(ku)) \to K(ku) \to K(KU)
$$
where $\eusm C^{KU}(ku)$ is the category of finite cell $ku$-module
spectra that become contractible when induced up to $KU$-modules.
Such spectra have finite Postnikov towers with layers that are induced
from finite cell $H\Z$-module spectra via the map $ku \to H\Z$, and
so it is reasonable to expect that a generalized form of the devissage
theorem in algebraic $K$-theory applies to identify $K(\eusm C^{KU}(ku))$
with $K(\Z)$.

\proclaim{Proposition 6.9}
If there is a fiber sequence $K(\Z) \to K(ku) \to K(KU)$ and $K(\Z)$
has telescopic complexity~$1$, at a prime $p\ge5$, then the algebraic
$K$-theory spectrum $K(KU)$ of the periodic topological $K$-theory
spectrum $KU$ has the same telescopic complexity as $K(ku)$, which is
$\ge 2$.
\endproclaim

\remark{Remark 6.10}
Unlike traditional elliptic cohomology, the spectrum $K(ku)$ is not
complex orientable.  For the map of $\SS$-algebras $ku \to H\Z$ is
$2$-connected, hence induces a $3$-connected map $K(ku) \to K(\Z)$;
cf.~\cite{W78, 1.1}, \cite{BM94, 10.9}.  So the composite map $\SS
\to K(ku) \to K(\Z)$ detects $\eta \in \pi_1(\SS)$.  This implies that
the unit map for $K(ku)$ cannot factor through the complex bordism
spectrum $MU$, since $\pi_1(MU) = 0$.  This should not be perceived as
a problem, however, as e.g.~also the topological modular forms spectrum $\tmf$
is not complex orientable.  It seems more likely that $K(\overline{KU})$
can be complex oriented, where $\overline{KU}$ is an ``algebraic closure''
of $KU$ in the category of commutative $\SS$-algebras.
\endremark


\enddocument